\hfuzz14pt
%\emergencystretch4mm
\input AHTOHFIE.STY
\def\AA{{\cal A}}
\def\BB{{\cal B}}
\def\CC{{\cal C}}
\def\DD{{\cal D}}
\def\FF{{\cal F}}

\def\PP{{\cal P}}
\def\RR{{\cal R}}
\def\NN{{\cal N}}
\def\KK{{\cal K}}
\def\LL{{\cal L}}
\def\MM{{\cal M}}
\def\QQ{{\cal Q}}
\def\SS{{\cal S}}

\def\itd#1#2{#1_1,\dots, #1_{#2}}

\let\wh\widehat
\def\cod{\hbox{\ssnine codim}\,}
%\def\Ends{{\rm Ends}\;}

%%%%%%%%%%%%%%%%%%%%%%%%%%%%%%%
\UDC{
512.567.5+%Структуры, связанные с универсальной алгеброй
512.544.44+%Характеристические и вполне характеристические подгруппы
512.543.22+%Тождественные соотношения в группах
512.544.33+%Нильпотентные и разрешимые группы
512.544.27+%Условия минимальности и максимальности
512.55+%Кольца и модули
519.171.1+%Общая теория графов
519.173.2%Топологические и метрические задачи теории графов. Планарность
}
\MSC{
06b05, %Lattices. Structure theory
08a35, %Algebraic structures. Automorphisms, endomorphisms
20e36, %Automorphisms of infinite groups
20f12, %Commutator calculus
20f16, %Solvable groups, supersolvable groups
20f18, %Nilpotent groups
20f19, %Generalizations of solvable and nilpotent groups
20e10, %Quasivarieties and varieties of groups
17a36, %Nonassociative rings and algebras. Automorphisms, derivations, other operators
17a30, %Nonassociative rings and algebras. Algebras satisfying other identities
05c10, %Planar graphs; geometric and topological aspects of graph theory
05c63  %Infinite graphs
}
\headline={\vbox{\hbox{\the\lhline}\hbox{\hss\the\rhline}\hss}}

\

\title
{Large and symmetric:\\
The Khukhro--Makarenko theorem on laws --- without laws
}

\author{Anton A. Klyachko\quad Maria V. Milentyeva}

\address{\myAddress\\ mariamil@yandex.ru}
\footnote{}{\hskip-0.5cm
The work of the first author was
supported by \RFBR11-01-00945.}

\medskip

\abstract{%
%\narrower
We prove a generalisation of the Khukhro--Makarenko theorem on
%том, что каждая группа, почти удовлетворяющая внешнему
%коммутаторному тождеству, содержит характеристическую подгруппу конечного
%индекса, удовлетворяющую этому тождеству.
large characteristic subgroups with laws.
This general fact implies new results on groups, algebras,
and even
graphs and other structures. Concerning groups, we obtain, e.g., a fact
in a sense
dual to the Khukhro--Makarenko theorem. A graph-theoretic corollary is an
analogue of this theorem in which planarity plays the role of a
multilinear identity. Also, we answer a question of Makarenko and
Shumyatsky.

%{\it Ключевые слова}:
%характеристические подгруппы, внешние коммутаторные тождества.
}
%%%%%%%%%%%%%%%%%%%%%%%%%%%%%%%%%%%%%%%%%%%%%%%%%%%%%%%%%%%%%%%%%%%%%%%%%%
%%%%%%%%%%%%%%%%%%%%%%%%%%%%%%%%%%%%%%%%%%%%%%%%%%%%%%%%%%%%%%%%%%%%%%%%%%

\s 0.
Introduction

The following short theorem generalised and strengthened various results
scattered through the literature (see, e.g., [BeK03], [BrNa04], and
Section 21.1.4 of [KaM82]).

\proclaim{Khukhro--Makarenko theorem \rm[KhM07a]}. If
a group has a finite-index subgroup satisfying an outer
commutator identity, then this group also has a characteristic
finite-index subgroup satisfying the same identity.

An \emph{outer \({\rm or} multilinear\) commutator identity} is an
identity of the form $[\dots[x_1,\dots,x_t]\dots]=1$ with some meaningful
arrangement of brackets, where all letters $x_1,\dots,x_t$ are different.
Examples of such identities are solvability, nilpotency,
centre-by-metabelianity, etc. A formal definition looks as follows. Let
$F(x_1,x_2,\dots)$ be a free group of countable rank.  An \emph{outer
commutator of weight 1} is just a letter $x_i$. An \emph{outer commutator
of weight $t>1$} is a word of the form
$w(x_1,\dots,x_t)=[u(x_1,\dots,x_r),v(x_{r+1},\dots,x_t)]$, where $u$ and
$v$ are outer commutators of weights $r$ and $t-r$, respectively. An
\emph{outer commutator identity} is an identity of the form $w=1$, where
$w$ is an outer commutator.

The Khukhro--Makarenko theorem has various applications (see, e.g.,
[KhM07b], [KhKMM09], [AST13] and references therein). Paper [KlM09]
contains a significantly simpler proof (than the original one) and a
better estimate of the index of the characteristic subgroup. At the same
time, in [KhM07b] and [KhM08], there were established some results (about
groups and algebras) similar to but not following from the
Khukhro--Makarenko theorem.

Paper [KhKMM09] was an attempt to clean up the mess; it contains a very
general proposition about groups with operators (in the sense of [Higg56],
see also [Kur62]) including as special cases all known and several new
results of this type. However, later, several facts that do not fit into
the general scheme of [KhKMM09] were discovered. Some of these facts are
rather delicate [MSh12], and others are quite simple (an impatient
reader may glance into the last section).
Generally, a Khukhro--Makarenko-like theorem looks as follows.

\proclaim{A theorem}.
If somewhere
there is something
{\small (in the classical case, a subgroup in a group)}
large
{\small (of finite index)}
and good
{\small (satisfying a multilinear identity)},
then
there is also something large, good, and symmetric
{\small (automorphism-invariant)}.

In this paper, we make another try to catch everything. In Section 1,
we prove the main theorem containing as special cases all known and
several new results similar to above (on groups,
algebras, graphs, and other objects).  The main idea is considering
``multilinear properties" instead of multilinear identities.

One of the new corollary of the main theorem
(Section 2) shows that, in the Khukhro--Makarenko
theorem, an identity, i.e. the triviality of a verbal subgroup, can be
replaced by something similar to the triviality. For
example, it is true that {\sl a group containing a finite-index subgroup
whose \the\year th derived subgroup is amenable (or periodic, or locally
finite, etc.)
contains a characteristic
finite-index subgroup with the same property}.

Another corollary (Section 3) can be considered as a result dual to the
Khukhro--Makarenko theorem. It shows, e.g., that,
{\sl for any
finite normal subgroup $N$ of any group $G$ of bounded
exponent, there exists a characteristic finite subgroup $H\nin G$ such
that the spectrum {\rm(i.e. the set of orders of all elements)} of the
quotient group $G/H$ is contained in the spectrum of $G/N$.  }

Section 4 contains similar results for algebras over fields.

Another proposition from Section 2 gives a positive answer to a question
of Makarenko and Shumyatsky and strengthens the main
theorem of [MSh12].

In Section 5, we show that some properties of graphs behave similarly to
multilinear commutator identities from the Khukhro--Makarenko theorem. For
example, it is true that, {\sl if some graph can be made planar by
removing a finite number of edges, then this finite set of edges can be
chosen invariant with respect to all automorphisms of the initial graph.}

As a reward to readers who reached the last section, we present 
two elementary-school problems in the subject.

The authors thank an anonymous referee and the editor for useful remarks.

%%%%%%%%%%%%%%%%%%%%%%%%%%%%%%%%%%%%%%%%%%%%%%%%%%%%%%%%%%%%%%%%%%%%%%%%%%
\s 1.
The main theorem

Recall that a \emph{semilattice} is a partially ordered set $\LL$ in
which any finite subset $\NN\subseteq\LL$ has a least upper
bound~${\sup\NN\in\LL}$. A \emph{directed semilattice} is a semilattice
which is a downward directed partially ordered set; this means that, for
any finite set $\NN\subseteq\LL$, there exists an element $\inf\NN\in\LL$
such that $\inf\NN\le N$ for any $N\in\NN$. Note that, in our notation,
$\sup$ is the least upper bound and $\inf$ is \emph{some} lower bound;
hopefully, this will not lead to a mess (see also the remark below, after
the definition of codimension). A semilattice $\LL$ is called
\emph{Noetherian} if all increasing chains in $\LL$ terminate, i.e. there
are no infinite chains of the form $N_1<N_2<\dots$, where $N_i\in\LL$. A
semilattice $\LL$ is called a \emph{lattice} if each finite subset
$\NN\subseteq\LL$ has a greatest lower bound
(which will be denoted $\inf\NN$).

We call a $t$-ary property (predicate) $\PP$ on
a semilattice $\LL$
\emph{(multi)monotone} if the property
$\PP(\itd N t)$, where $N_i\in\LL$, implies
$\PP(\itd {N'} t)$
for any $N'_i\leq N_i$.
We call the property $\PP$ \emph{multilinear} if, for any $i$,
\newline
$\PP(\itd N {i-1}, N'_i,N_{i+1},\dots,N_t)$ and
$\PP(\itd N {i-1}, N''_i,N_{i+1},\dots,N_t)$
imply
\newline
$\PP(\itd N {i-1},\sup(N'_i,N''_i),N_{i+1},\dots,N_t)$.

We need also dual notions. A predicate $\PP$ on
a semilattice $\LL$ is called
\emph{comonotone} if
$\PP(\itd N t)$, where $N_i\in\LL$, implies
$\PP(\itd {N'} t)$
for any $N'_i\geq N_i$.
We call a predicate $\PP$ \emph{comultilinear} if, for any $i$,
\newline
$\PP(\itd N {i-1}, N'_i,N_{i+1},\dots,N_t)$ and
$\PP(\itd N {i-1}, N''_i,N_{i+1},\dots,N_t)$
imply
\newline
$\PP(\itd N {i-1},\inf(N'_i,N''_i),N_{i+1},\dots,N_t)$
for some lower bound $\inf(N'_i,N''_i)$ (we shall use the word 
\emph{colinear} if $t=1$).

We say that an endomorphism semigroup $\Phi\subseteq\End\LL$ of a
semilattice $\LL$ \emph{preserves a property $\PP$} (or the property~$\PP$
is \emph{$\Phi$-invariant}) if $\PP(\itd N t)$ implies
$\PP(\phi(N_1),\dots,\phi(N_t))$ for any $N_i\in\LL$ and $\phi\in\Phi$.
For example, the property $\RR(X,Y,Z)=\bigg(X=\sup(Y,Z)\bigg)$ is
$(\End\LL)$-invariant by the definition of \emph{endomorphisms of a
semilattice}.  An element~$N$ of a semilattice is called
\emph{$\Phi$-invariant} if $\phi(N)\le N$ for all $\phi\in\Phi$ (this, in
particular, means that $\phi(N)=N$ for all $\phi\in\Phi$ if the semigroup
$\Phi$ is a subgroup of $\Aut\LL$).

The following assertion is a natural generalisation of the lemma from
[KlM09] (which is about the lattice of normal subgroups),
see also Lemma 1 from [KhKMM09]
(which is about the lattice of normal subgroups in a multi-operator group).

\Lemma 1. %\label{ta}
Suppose that $\MM$ is a directed semilattice with a largest
element $\sup\MM$,\ \ $\PP$ is a monotone multilinear
$t$-ary predicate on $\MM$,\ \ $m$ is a positive integer, and
$\NN\subseteq \MM$ is a finite subset of $\MM$ such that
$$
\PP(\underbrace{N,N,\dots,N}_{m\rm\;times},
\sup\MM,\sup\MM,\dots,\sup\MM)
\quad\hbox{is true for all $N\in\NN$}.
$$
Then
$$ %\label{ea}
\PP(\underbrace{\widehat N,\wh N,\dots,\wh N}_{m-1\rm\;times},
\wh G,\wh G,\dots,\wh G),
\quad\hbox{where } \wh N=\inf\NN \hbox{ and }
\wh G= \sup\NN, \hbox{ is true too}.
$$

\Proof
Since $\PP$ is monotone, $\wh G\leq\sup\MM$, and $\wh N\leq N$ for
all $N\in\NN$, we have
$$
\PP(\underbrace{\widehat N,\wh N,\dots,\wh N}_{m-1\rm\;times},
N,\wh G,\dots,\wh G)
\quad\hbox{is true for all $N\in\NN$}.
$$
Now the multilinearity (to be more precise, the linearity with respect to
the $m$th argument) implies the assertion of Lemma 1.

\medskip

Let $\LL$ be a Noetherian directed
semilattice and let $\Phi\subseteq\End \LL$
be a semigroup of endomorphisms of~$\LL$.
A function $\cod\!\:\LL \to \R$ is called a \emph{(generalised)
$\Phi$-codimension} if it has the following properties:
\item{1)}
$\cod N_1 \le\cod N_2$ if $N_1 \geq N_2$;
\item{2)}
$\cod \phi(N)\leq \cod N$ for any $N\in\LL$ and $\phi\in \Phi$;
\item{3)}
$\cod \inf(N_1,N_2)\leq \cod N_1+\cod N_2$ for any
$N_1,N_2\in \LL$ and some lower bound $\inf(N_1,N_2)$;
\item{4)}
in any family $\cal N \subseteq \LL$, there exists
$r\leq \max\limits_{N\in \cal N}\cod N+1$ elements $\itd N r$ such
that
$
\sup{\cal N}=\sup(\itd N r).
$

\enditem
This definition of codimension is a natural generalisation of the
corresponding notion from [KhKMM09]
(which is about the lattice of normal subgroups of a multi-operator
group). When we are talking about a codimension on a semilattice
we suppose that the symbol $\inf$ always denote a lower bound satisfying
condition 3).

\proclaim{Main theorem}. %\label{tb}
Suppose that $\LL$ is a Noetherian directed
semilattice, $\Phi \subseteq \End \LL$ is a semigroup of its
endomorphisms, and $\PP$ is a multi-monotone multilinear $t$-ary
$\Phi$-invariant predicate on $\LL$.
Then, if there exists an element
$N\in\LL$ with the property $\PP(N,\dots, N)$,
then there exists an element $H\in \LL$ such that
\item{\rm 1)}
$H$ has the same property: $\PP(H,\dots,H)$;
\item{\rm 2)}
$H$ is $\Phi$-invariant;
\item{\rm 3)}
if $\phi(N)\le J$ for any $\phi\in \Phi$
and some $J\in\LL$, then $H\le J$;
\item{\rm 4)}
if $\LL$ is a lattice (i.e. each finite set has a greatest
lower bound)
and $\Phi$ consists of lattice endomorphisms (i.e. the mappings
commuting with the taking the greatest lower bounds of finite sets),
then $H$ is contained in the sublattice generated by the set
$\{\phi(N)\;;\;\phi\in\Phi\}$;
\item{\rm 5)}
if $\cod\!\:\LL \to \R$ is a generalised $\Phi$-codimension, then
$
\cod H \leq f^{t-1}(\cod N),
$
where $f^k(x)$ is the $k$th iteration of the function $f(x)=x(x+1)$.

\Proof
Since $\LL$ is Noetherian, it contains an element
$G_1=\sup\limits_{\phi\in \Phi}\phi(N)$
and
$$
G_1=\sup(\phi'_0(N),\dots,\phi'_{p_1}(N))
\qbox{for some endomorphisms $\phi_i' \in \Phi$}.
$$
Note that $G_1$ is $\Phi$-invariant:
$
\phi(G_1)=
\sup(\phi\phi'_0(N),\dots,\phi\phi'_{p_1}(N))
\le\sup\limits_{\phi\in\Phi}\phi(N)=G_1.
$
Besides, $G_1\le J$ (if $J$ is such as in assertion 3) of the
theorem) and, for any codimension function on $\LL$, we have
$$
p_1\leq l_0\:=\cod N
\qbox{by property~4) of codimension}.
$$
Put $N_1=\inf(\phi'_0(N),\dots,\phi'_{p_1}(N))$.
By properties~2) and~3) of $\cod$, we have
$$
l_1\:=\cod N_1 \leq (p_1+1)\cod N =(p_1+1)l_0\leq(l_0+1)l_0=f(l_0).
$$
According to Lemma 1 (applied to the semilattice $\MM=\LL$)
we have the property
$$
\PP(N_1,\dots,N_1,G_1).
$$
Similarly, we can choose elements
$$
G_2=\sup_{\phi\in\Phi}\phi(N_1)=\sup(\phi''_0(N_1),\dots,\phi''_{p_2}(N_1))
\ \mbox{  and  }\
N_2=\inf(\phi''_0(N_1),\dots,\phi''_{p_2}(N_1)),
\qbox{where
$\phi''_i\in\Phi$}.
$$
Clearly,
$$
\eqalign{
&G_2\le G_1\le J
\qbox{(because $N_1\le\phi'_0(N)$
and, hence,
$G_2=\sup\limits_{\phi\in\Phi}\phi(N_1)\le
\sup\limits_{\phi\in\Phi}\phi\phi'_0(N)\le
\sup\limits_{\phi\in\Phi}\phi(N)=G_1$)}
\qqbox{and}
\cr
&p_2\leq \cod N_1=l_1\leq f(l_0)
\qbox{(by properties 2) and 4) of $\cod$)}.
}
$$
The element $G_2$ is obviously $\Phi$-invariant
(for the same reasons as $G_1$)
and
we have the estimates
$$
\cod G_2 \le
\cod\phi_0''N_1\le
\cod N_1= l_1\le f(l_0)
\qqbox{and}
l_2\:=\cod N_2\leq (p_2+1)\cod N_1=(p_2+1)l_1\leq f(l_1)\leq f(f(l_0)).
$$
Again by Lemma 1 (applied to the semilattice
$\MM=\{X\in\LL\;|\; X\le G_1\}$),
we obtain
$$
\PP(N_2, \dots, N_2,G_2,G_2).
$$
Continuing in the same manner, at the $t$th step,
we obtain a $\Phi$-invariant element
$$
G_t=\sup_{\phi\in \Phi}\phi(N_{t-1})=
\sup(\phi^{(t)}_0(N_{t-1}),\dots,\phi^{(t)}_{p_t}(N_{t-1}))
\hbox{ and an element }
N_t=\inf(\phi^{(t)}_0(N_{t-1}),\dots,\phi^{(t)}_{p_t}(N_{t-1})),
$$
where $\phi^{(t)}_i\in \Phi$.
We have the required property $\PP(G_t,\dots, G_t)$ and the inequalities
$$
G_t\le J,
\qquad
\cod G_t\leq \cod N_{t-1}=l_{t-1}\leq f(l_{t-2})\leq f(f(l_{t-3}))
\leq \dots \leq f^{t-1}(l_0).
$$
Thus, the element $H=G_t$ is as required and the theorem is proven
(Assertion 4 is obviously satisfied by the construction).

\medskip

The following lemma makes it possible to construct new multilinear
properties from known ones.

\proclaim{Composition Lemma}.
Suppose that, on a lattice, there is a multilinear monotone predicate
%\newline
$\QQ(M_1,\dots,M_k)$ and a tuple of predicates
$
\RR=\left\{\RR_i
\pmatrix{
X_1,\dots,X_{l}\cr
Y
}\right\}
$
which are multilinear and monotone with respect to the first row
(i.e. for any given second row)
and colinear and comonotone with respect to the second row
(i.e. for any given first row).
Then the predicate
$$
\QQ\o\RR(N_1,\dots,N_{kl})=
\(\exists M_1,\dots,M_k\quad
\QQ(M_1,\dots,M_k) \hbox{ and, for $i\in\{1,\dots,k\}$, }\
\RR_i
\pmatrix{
N_{(i-1)l+1},\dots,N_{li}\cr
M_i
}
\)
$$
called the \emph{composition} of the predicates $\QQ$ and $\RR$ is
multilinear and monotone.

\Proof
Let us verify
that the composition is
monotone, e.g., with respect to the first argument. If
$N_1'\le N_1$ and the property ${\QQ\o\RR(N_1,\dots,N_{kl})}$ holds
(for some $M_1,...,M_k$ from the definition of composition), then
${\QQ\o\RR(N_1',\dots,N_{kl})}$ holds also (with
the very same $M_i$) because $\RR_1$ is monotone with respect
to the first row.

Let us verify the multilinearity of $\QQ\o\RR$,
e.g., with respect to the first
argument. Suppose that the properties $\QQ\o\RR(N_1',N_2,\dots,N_{kl})$
and
$\QQ\o\RR(N_1'',N_2,\dots,N_{kl})$ hold, i.e. we have
$$
\eqalign{
&\RR_1
\pmatrix{
N_1',N_2,\dots,N_l\cr
M_1'
},\
\RR_2
\pmatrix{
N_{l+1},\dots,N_{2l}\cr
M_2'
},
\dots,
\qqbox{\phantom{and}}
\QQ(M_1',\dots,M_k'),
\cr
&\RR_1
\pmatrix{
N_1'',N_2,\dots,N_l\cr
M_1''
},\
\RR_2
\pmatrix{
N_{l+1},\dots,N_{2l}\cr
M_2''
},
\dots\phantom,
\qqbox{and}
\QQ(M_1'',\dots,M_k'')
}
$$
for some $M_1',\dots,M_k',M_1'',\dots,M_k''$.
We claim that then we have the properties
$$
\eqalign{
&\RR_1
\pmatrix{
\sup(N_1',N_1''),N_2,\dots,N_l\cr
\sup(M_1',M_1'')
},\
\RR_2
\pmatrix{
N_{l+1},\dots,N_{2l}\cr
\inf(M_2',M_2'')
},
\dots,
\qqbox{ and}
\QQ(\sup(M_1',M_1''),\inf(M_2',M_2''),\dots,\inf(M_k',M_k''))
}
$$
(i.e. the property $\QQ\o\RR(\sup(N_1',N_1''),N_2,\dots,N_k)$ holds
with
\newline
$\sup(M_1',M_1'')$, $\inf(M_2',M_2'')$,\dots, $\inf(M_k',M_k'')$
in the roles of $M_1,\dots,M_k$,
respectively). Indeed,
\-
the first property ($\RR_1(\dots)$) holds because
$\RR_1$ is linear with respect to the first element of the first
row and comonotone with respect to the second row;
\-
the second property ($\RR_2(\dots)$) holds by virtue of the colinearity
of $\RR_2$ with respect to the second row;
\-
\dots
\-
the last property
($\QQ(\dots)$) holds because $\QQ$ is linear with respect to
the first argument and monotone with respect to the other arguments.

\enditem
The lemma is proven.

\medskip

The rest of the paper deals with applications of the main theorem to
groups, algebras, graphs, and other structures.
In all applications, the semigroup $\Phi$ is taken to be a natural
subgroup of $\Aut\LL$.

%%%%%%%%%%%%%%%%%%%%%%%%%%%%%%%%%%%%%%%%%%%%%%%%%%%%%%%%%%%%%%%%%%%%%%%%%%
\s 2.
The lattice of large normal subgroups

Recall that an abstract class of groups $\KK$ is called
\emph{radical} (or \emph{Fitting}) if it is closed
with respect to normal
subgroups and finite products of normal subgroups, i.e.
\item {1)}
any normal subgroup of a group from $\KK$ lies in $\KK$;
\item {2)}
a group decomposable
into a product of two normal subgroups lying in $\KK$ belongs to
$\KK$.

\enditem
A \emph{coradical class} (or a \emph{formation}) is an abstract
class of groups $\KK$ closed with respect to homomorphic images and
subdirect products, i.e. such that:
\item {$1'$)}
any quotient group of a group
from $\KK$ lies in $\KK$;
\item {$2'$)}
any subdirect product of two groups
lying in $\KK$ also lies in $\KK$.
\enditem
More details about radical
and coradical classes can be found, e.g., in book~[She78].

The following classes of groups are
\emph{radical formations}, i.e. they
are both radical and coradical:
\- finite groups;
\- finite $p$-groups;
\- locally finite groups
(radicality follows from a theorem of O. Yu. Schmidt:
{\sl an extension of a locally finite group by a locally finite
is locally finite itself}, see [KaM82]);
\- periodic groups;
\- Noetherian groups;
%\- группы, нётеровы по нормальным подгруппам;
\- Artinian groups;
\- nilpotent groups
(radicality follows from the Fitting theorem, see [KaM82]);
\- locally nilpotent groups
(radicality follows from Plotkin's theorem, see [KaM82]);
\- solvable groups;
\- almost solvable groups;
\- locally polycyclic groups
(radicality follows from Theorem 18.1.2 of [KaM82]);
\- groups satisfying nontrivial identities;
\- groups without nonabelian free subgroups;
\- amenable (discrete) groups;
\- \dots
\enditem
Other examples of coradical classes are all varieties of groups,
the class of all binary finite groups, groups with the
maximality (or minimality) conditions for normal
subgroups, etc.

%Примером радикального класса, не являющегося корадикальным, может служить
%класс всех групп без кручения.
% замкнутость относительно произведений: элемент кон. порядка
% в AB порождает подгруппу, трив. пересекающуюся с A,... нет,!!!

\proclaim{Large-subgroup theorem}.
Let $N$ be a normal subgroup of a group $G$ such that
the quotient group $G/N$
satisfies the maximality condition for
normal subgroups.
Then $G$ contains characteristic
%и даже инвариантные относительно всех сюръективных эндоморфизмов
subgroups
$H_1, H_2,\dots$ such that
\item{\rm 1)}
the quotient groups $G/H_t$ lie in the coradical class (formation) $\FF$
generated by $G/N$;
moreover, the subgroups~$H_t$ belongs to the lattice (of subgroups
of $G$) generated by the images of $N$ under all
automorphisms of $G$;
%если $\Phi\subseteq\Aut G$;
\item{\rm 2)}
for any multilinear commutator word $w$ of degree
$k\le t$,
the group
$w(H_t,\dots,H_t)$ is contained in the radical class~$\RR_w$
generated by the group $w(N,\dots,N)$;
\item{\rm 3)}
if $\cod$ is a generalised codimension on
the lattice of subgroups such that the corresponding quotients
lie in $\FF$, then $\cod H_t \leq f^{t-1}(\cod N)$, where $f^k(x)$
is the $k$th iteration of the function ${f(x)=x(x+1)}$.

\Proof
Let $\itd K t $ be normal subgroups of $G$ and let $w(\itd x t)$ be an
outer commutator. Then

\item{a)}
the subgroup $w(\itd K t)\:=\gp{w(\itd h t) : h_i\in K_i}$
is normal in $G$;
\item{b)}
$w(\itd K t)=[u(\itd K r), v(K_{r+1},\dots,K_t)]$
if $w(\itd x t)=[u(\itd x r), v(x_{r+1},\dots,x_t)]$;
\item {c)}
$w(\itd K {i-1}, \prod\limits_{N\in \cal N}N,K_{i+1}, \dots, K_t)=
\prod\limits_{N\in \cal N}w(\itd K {i-1}, N,K_{i+1}, \dots, K_t)$
for any family~$\NN$ of normal subgroups of $G$.

\enditem
These facts are well known and easy to prove by induction.

Let us apply the main theorem to
the lattice $\LL$ (of normal subgroups of $G$)
generated by the images of $N$ under all
automorphisms of $G$.
%факторгруппы по которым
%лежат в $\FF$
(This lattice is Noetherian and even the entire formation $\FF$ consists of
groups which are
Noetherian with respect to normal subgroups.)
Put $\Phi=\Aut G$,
%моноид сюръективных эндоморфизмов
%группы $G$ (такие эндоморфизмы индуцируют эндоморфизмы
%полурешётки $\LL$),
and let $\PP(\itd N t)$ be the following property:

{\narrower\narrower\narrower\narrower\narrower\narrower\narrower

\noindent\sl
for each multilinear commutator word $w$ of degree at most $t$, the
group $w(\itd N k)$
belongs to~$\RR_w$.

}
\noindent
The predicate $\PP$ is monotone because
radical classes are closed with respect to normal
subgroups; the multilinearity follows from the closedness of radical
classes with respect to products of normal subgroups and
property~c) of outer commutators. The theorem is proven.

\medskip
\noindent
The large-subgroup theorem generalises the Khukhro--Makarenko theorem in
three directions:

\-
the finiteness
of the quotient groups $G/N$ and $G/H$
is replaced by their belonging
to any given formation with maximality condition for
normal subgroups and the index (to be more precise, the logarithm of
index) is replaced by an abstract codimension; this
generalisation is not new; it was obtained in [KhKMM09] for the first time;

\-
the identity (i.e. the triviality of a verbal subgroup) is replaced by
the
belonging of this verbal subgroup to any radical
class; e.g., our theorem shows that {\sl a group containing
a finite-index subgroup whose 100th derived subgroup is periodic
contains a characteristic finite-index subgroup with the same property};

\-
finally, instead of one multilinear word $w$, we consider all
multilinear words at once; this gives a substantial gain in estimation
if, e.g., we want to
construct a characteristic finite-index subgroup
satisfying \emph{all} multilinear identities of degree
at most 100  satisfied by a given finite-index subgroup (of course,
such a characteristic subgroup
can be constructed by iterated applications of the Khukhro--Makarenko
theorem but this results in a very bad estimate of index).

\enditem
The following simple facts make it possible to apply
the composition lemma, generalise the Khukhro--Makarenko theorem in
yet another direction, and answer a question of Makarenko and
Shumyatsky.

\proclaim{Quotient-group lemma}.
The following two properties $\AA(N,M)$ and $\BB(N,M)$ of
pairs of normal subgroups of a group $G$
can
be written in the form $\RR\pmatrix{N,\dots,N\cr M}$, where the predicate
$\RR$ on the lattice of normal subgroups is
monotone and multilinear with
respect to the first row and
comonotone and colinear with respect to the second row:
$$
\eqalign{
&\AA(N,M)=\bigg(N/(N\cap M) \hbox{ satisfies a (given) outer
commutator identity $w=1$}\bigg),
\cr
&\BB(N,M)=\bigg(N/(N\cap M) \hbox{ belongs to a (given)
radical formations $\FF$}\bigg).
\cr
}
$$

\Proof
For the property $\AA$, the predicate
$$
\RR\pmatrix{N_1,\dots,N_t\cr M}=
\bigg(w(N_1,\dots,N_t)\subseteq M\bigg)
$$
obviously satisfies all conditions.
As for $\BB$, this property itself
can play the role of $\RR$:
$$
\RR\pmatrix{N_1\cr M}=\BB(N_1,M).
$$
This predicate is linear and monotone with respect to the first
row because of the radicality of $\FF$; colinearity and
comonotonicity with respect to the second row follow from the coradicality
of $\FF$.

\proclaim{Extension lemma}.
Suppose that $\QQ(M_1,\dots,M_l)$ is a monotone multilinear predicate
on the lattice of normal subgroups of a group $G$,\ \
$w$ is an outer commutator word of degree $d$,
and $\FF$ is a radical formation.
Then the following two properties $\CC(N)$ and $\DD(N)$ of
normal subgroups
of $G$ can
be written in the form $\PP(N,\dots,N)$, where the predicate
$\PP(N_1,\dots,N_t)$ on the lattice of normal subgroups is
monotone and multilinear, $t=ld$ for the property
$\CC$, and $t=l$ for the property $\DD$:
$$
\eqalign{
&\CC(N)=
\bigg(
\exists M\nin G\quad \ M\subseteq N,\  \QQ(M,\dots,M),\
\hbox{ and }
N/M \hbox{ satisfies the identity $w=1$}
\bigg),
\cr
&\DD(N)=
\bigg(
\exists M\nin G\quad \ M\subseteq N,\  \QQ(M,\dots,M),\
\hbox{ and }
N/M \in\FF
\bigg).
\cr
}
$$

\Proof
Clearly,
the property $\CC(N)$ can be rewritten in the form
$$
\CC(N)=\big(\exists M\nin G\quad \AA(N,M)\hbox{ and } \QQ(M,\dots,M)\big),
\qbox{ where $\AA$ is from the quotient-group lemma}.
$$
Note that the last formula is equivalent to
$$
\CC(N)=
\bigg(
\exists M_1,\dots,M_l\nin G\quad
\AA(N,M_1),\ \dots,\
\AA(N,M_l),
\hbox{ and }
\QQ(M_1,\dots,M_l)
\bigg).
$$
To verify this (in the non-obvious direction), it suffice to put
$M=\bigcap M_i$ and recall that the variety of groups satisfying the
identity
$w=1$ is closed with respect to subdirect products and the property $\QQ$
is monotone.
By the quotient group lemma, we can rewrite $\CC$ in the form
$$
\CC(N)=
\bigg(
\exists M_1,\dots,M_l\nin G\quad
\RR\pmatrix{N,\dots,N\cr M_1},\ \dots,\
\RR\pmatrix{N,\dots,N\cr M_l},\quad
\hbox{and }
\QQ(M_1,\dots,M_l)
\bigg),
$$
where the predicate
$\RR$ on the lattice of normal subgroups is
monotone and multilinear with respect to the first row,
and
comonotone and colinear with respect to the second row.

Application of the composition lemma completes the proof.
Similar arguments prove the assertion about
$\DD$.

\medskip

The following theorem was proven in [MSh12] for the case, where the group
$G$ is locally finite and each class $\KK_i$ is the class of all
locally nilpotent groups.

\proclaim{Series theorem}.
Suppose that a group $G$ contains a finite-index subgroup $N$ having
normal (in $G$) series
$$
\1=A_0\subseteq\dots\subseteq A_n=N
$$
such that each quotient $A_i/A_{i-1}$ either satisfies a multilinear
commutator identity $w_i=1$ of weight $t_i$ or lies in a radical 
class~$\KK_i$; moreover,
all these classes, except possibly $\KK_1$, are also coradical. Then $G$
contains a characteristic  subgroup $H$ with the same property
(i.e. with a series of the same length with quotients satisfying the same
identities or lying in the same classes) and 
$\log_2|G:H| \leq f^{t-1}(\log_2|G:N|)$, where $f^k(x)$
is the $k$th iteration of the function 
${f(x)=x(x+1)}$ and $t=\prod t_i$.

\Proof
The extension lemma and an obvious induction show that the presence
of such a normal series
in a normal subgroup $N$ can be written in the form
$\PP(N,\dots,N)$, where $\PP$ is a multilinear monotone $t$-variable 
predicate on the lattice of normal subgroups of $G$.
It remains to apply the main theorem.

\medskip

Note that, in the case where the group
$G$ is locally finite and each class $\KK_i$ is the class of all locally
nilpotent groups, in
[MSh12] a stronger proposition was proven:  the group
$G$ has a characteristic finite-index subgroup having a
\emph{characteristic} series with the required property.  In the general
case, such a strengthening is impossible as the following example due
to Yves Cornulier shows.

\Example {\rm [Corn13]}.
{\sl There exists a group with a normal abelian countable-index
subgroup but without characteristic abelian countable-index
subgroups.}
\newline
Take a countable-dimensional vector space $V$ with a basis $\{e_q\}$,
where $q\in\Q$, over a finite field $K$ and consider the group~$G$
of  ``unitriangular" operators, i.e. such operators $g$
in $V$ that $ge_q-e_q\in\gp{\{e_r:r<q\}}$. Now, consider the
subgroup $H\subset G$ consisting of matrices $A$ with the following
property: for any real $r$, there is only finitely many nonzero entries
$a_{pq}$ such that $p\ne q$ and either $p>r$ or $q<r$.
This group $H$ has an abelian normal countable-index subgroup
consisting of matrices $A$ whose all nonzero nondiagonal entries
$a_{pq}$ are such that either $p<0$ or $q>0$. However, $H$ has no
nontrivial characteristic abelian subgroups.  Indeed, if
$1\ne h\in N\nin H$, then the commutator $c$ of $h$ and any transvection
lies in $N$ and has only finitely many nonzero nondiagonal entries; hence,
$c$ belongs to a finite-dimensional unitriangular group
$\UT_n(K)$ which is nilpotent and, therefore, any its nontrivial normal
subgroup nontrivially intersect the centre (which consists of
transvections).
Thus, any nontrivial normal subgroup of $H$ contains a transvection.
It remains to note that the automorphism group of $H$ acts transitively
on the transvections, i.e. any nontrivial characteristic subgroup of $H$
must contain all transvections and, hence, cannot be abelian.

%%%%%%%%%%%%%%%%%%%%%%%%%%%%%%%%%%%%%%%%%%%%%%%%%%%%%%%%%%%%%%%%%%%%%%%%%%
\s 3.
The lattice of small subgroups. A dual of the Khukhro--Makarenko theorem

Consider a universal positive closed first-order formula
in the group language, e.g.,
$$
(\forall x)(\forall y)\
\Big(\(x^3=y^3\ \wedge\ (xy)^4=(yx)^4\)\ \ \vee\ \ (xy)^{\the\year}=1\
\ \vee\ \ [x,y]^5=1\Big).
$$
Such a formula defines a class of groups consisting of groups where the
formula holds. For instance, the formula
${(\forall x)\;(x^2=1\ \vee\ x^3=1)}$ holds in the symmetric group of
order six, but does not hold in the abelian group of order six.

The following theorem can be considered as dual to
the Khukhro--Makarenko theorem.

\proclaim{Finite-subgroup theorem}.
If a group $G$ has a finite normal subgroup such that, in the quotient
group,
a given universal positive closed first-order formula holds,
then $G$ has a characteristic finite subgroup with the same
property.

We shall prove a more general assertion similar to
the large-subgroup theorem.
Consider the following
property~$\DD(N)$ of a normal subgroup $N$ of a group $G$:
$$
(\forall x)(\forall y)\dots
\(
\SS(x,y,\dots)
\imp
\bigvee_{i=1}^t
(G/N,xN,yN,\dots)\in\FF_i
\),
\eqno{(*)}
$$
where
$\SS$ is some $(\Aut G)$-invariant property of a tuple of elements
of $G$ and $\FF_i$ are some formations of groups with marked
elements.

A \emph{formation of groups with marked elements} is a
class of tuples $(H,h_1,h_2,\dots)$,
where $H$ is a group and $h_i\in H$,
closed with respect to homomorphic images and subdirect products (which
are defined naturally). For example, the class
of amenable groups with two marked elements whose commutator lie
in the centre is a formation.

Properties (of normal subgroups) of the form $(*)$ are called
\emph{$t$-disjunctive}.

\proclaim{Small-subgroup theorem}.
If all groups from
the radical class generated by a normal subgroup $N$ of a group $G$
satisfy the minimality conditions for
normal subgroups, then $G$ contains
characteristic subgroups $H_1, H_2,\dots$ such that
\item{\rm 1)}
$H_t$ lie in the radical class $\RR$,
generated by
the group $N$;
moreover, $H_t$ are contained in the lattice (of subgroups
of~$G$) generated by all automorphic images of $N$;
\item{\rm 2)}
the quotient group $G/H_t$ satisfies all $t$-disjunctive properties
satisfied by $G/N$;
\item{\rm 3)}
if, in addition, $\cod$ is a generalised codimension
{\rm(which may be called a dimension in this case)},
on the
lattice dual to the lattice of subgroups lying in $\RR$, then
$
\cod H_t \leq f^{t-1}(\cod N),
$
where $f^k(x)$ is the $k$th iteration of the function ${f(x)=x(x+1)}$.

\Proof
Each $t$-disjunctive property $\DD(N)$ can be rewritten in the form
$\DD(N)=\PP_\DD(N,\dots,N)$, where $\PP_\DD(N_1,\dots,N_t)$ is
the following property of a tuple of normal subgroups:
$$
(\forall x)(\forall y)\dots
\(
\SS(x,y,\dots)
\imp
\bigvee_{i=1}^t
(G/N_i,xN_i,yN_i,\dots)\in\FF_i
\).
$$

Now, it suffice to apply the main theorem to the
lattice
$\LL$
(of subgroups of~$G$) generated by
all images of $N$ under
automorphisms of $G$; but the order on $\LL$
is
opposite to the natural one:
$
A\le B \hbox{ if } A\supseteq B.
$
This lattice is Noetherian, because it consists of normal subgroups
that are Artinian with respect to normal subgroups. The
semigroup $\Phi$ is the automorphism group of $G$ in this case
and the predicate $\PP$ is the conjunction of all predicates $\PP_\DD$,
where $\DD$ runs over all $t$-disjunctive properties satisfied by $N$.

The predicate $\PP$ is monotone because
formations are closed with respect to quotient groups,
the multilinearity follows from the closedness of formations
with respect to subdirect products.
Let us verify the linearity, e.g., with respect
to the first argument. We have to
prove that the property $\PP_\DD(N_1'\cap N_1'',N_2,\dots)$ holds
whenever
the properties $\PP_\DD(N_1',N_2,\dots)$ and
$\PP_\DD(N_1'',N_2,\dots)$
hold. Thus, we known that
for any set of elements $g,h,\dots\in G$ with property $\SS$,
either
\-
$(G/N_i,\,gN_i,\,hN_i,\dots)\in\FF_i$ for some $i\ge2$
\-
or the formation $\FF_1$ contains
two group with marked elements:
\newline
$(G/N_1',\,gN_1',\,hN_1',\dots)$ and $(G/N_1'',\,gN_1'',\,hN_1'',\dots)$
\newline
and, hence, it contains
their subdirect product
$(G/(N_1'\cap N_1''),\,g(N_1'\cap N_1''),\,h(N_1'\cap N_1''),\dots)$.

\enditem
This means that the property $\PP_\DD(N_1'\cap N_1'',N_2,\dots)$ holds and
the theorem is proven.

\medskip
The finite-subgroup theorem is obtained from the small-subgroup theorem
by taking all formations $\FF_i$ to be the class of all
groups whose marked elements satisfy a system of
equations (depending on $i$). The following is a special case
of the finite-subgroup theorem.

\proclaim{Spectrum theorem}.
For any finite normal subgroup $N$
of a
group $G$ of bounded exponent,
there exists a characteristic finite subgroup $H$ such that
the
spectrum {\rm(i.e. the set of orders of all elements)}
of $G/H$
is contained in the spectrum of $G/N$.

\Proof
It suffices to apply the finite-subgroup theorem to the formula
$
\
\forall x
\
\bigvee\limits_{i=1}^t
\(x^{n_i}=1\),
%\quad
$
{where $\{n_1,\dots,n_t\}$ is the spectrum of $G/N$.}

\medskip
The order of the characteristic subgroup $H$ can be explicitly
estimated via the order of $N$ and the cardinality of the spectrum
of~$G/N$ (because logarithm of the order of a subgroup is a natural
example of codimension on the lattice of finite normal
subgroups). The word ``finite" (the both occurrences) in the spectrum
theorem can be replaced by ``Artinian" or, e.g., ``Chernikov", etc.  The
word ``spectrum" (the both occurrences) can be replaced, e.g., by
``spectrum of the commutator subgroup"; to show this we may
add to the property
$\SS(x)$ the condition that $x$ lies in the commutator subgroup.

In conclusion, we give an example showing that
the finite-subgroup theorem cannot be extended to arbitrary
(non-universal) positive first-order formulae.

\Example.
The group $G=\gp a_2\times B$, where $B$ is an abelian divisible group
with infinitely many elements of order two
(e.g., $B=\(\Z_{2^\infty}\)^\infty$), has an obvious finite
normal
subgroup $N=\gp a_2$ such that, in the quotient group, all
elements are squares (i.e., the formula
$\forall x\ \exists y\ x=y^2$ holds). But there is no characteristic
finite subgroup with the same property.
Indeed, such a
characteristic subgroup $H$ cannot lie in $B$, obviously. Take an element
$(a,b)\in H$ and consider its images under automorphisms that
fix elements of $B$ and map $a$ into $(a,x)$, where
$x$ is an element of order two from~$B$. These images $(a,bx)$ form
an infinite subset of $H$.

%%%%%%%%%%%%%%%%%%%%%%%%%%%%%%%%%%%%%%%%%%%%%%%%%%%%%%%%%%%%%%%%%%%%%%%%%%
\s 4.
The lattices of ideals and subspaces

The term \emph{algebra} in this section means not necessarily
associative algebra over a field. A \emph{characteristic} subspace of an
algebra is a subspace invariant with respect to all automorphisms of this
algebra.

\proclaim{Large-subspace theorem}.
Let $N$ be a subspace of an algebra $G$ such that
either
\-
$N$ is of finite codimension,
\-
$N$ is a left ideal and the
 $G$-module $G/N$ is Noetherian,
\-
or $N$ is a two-sided ideal and
the quotient algebra from
$G/N$ satisfies
the maximality condition for two-sided ideals.
\enditem
Then $G$ contains a characteristic subspaces
$H_1, H_2,\dots$ such that
\item{\rm 1)}
$H_t$ belong to the lattice (of subspaces
of $G$) generated by the images of $N$ under all
automorphisms of $G$;
in particular, the subspaces $H_t$ are ideals
(one-sided or two-sided) if $N$ is an ideal,
the quotient algebras (quotient modules)
$G/H_t$ lie in the formation $\FF$,
generated by the algebra (module) $G/N$;
%если $\Phi\subseteq\Aut G$;
\item{\rm 2)}
for any multilinear element
$w(x_1,\dots,x_n)$ of the free (nonassociative) algebra of rank
$n\le t$,
the set
$w(H_t,\dots,H_t)$ is contained in the linear hull of a
finite number of images of the
set $w(N,\dots,N)$ under automorphisms of $G$;
\item{\rm 3)}
If, in addition, $\cod$ is
either
the usual codimension (of a subspace of $G$)
or
a generalised codimension on the
lattice of ideals such that the
corresponding
quotient algebras (modules) lie in $\FF$
then $\cod H_t \leq f^{t-1}(\cod N)$, where $f^k(x)$ is the $k$th
iteration of the function ${f(x)=x(x+1)}$.

\noindent{\bf The proof}
almost literally repeats the proof of the large-subgroup theorem;
only obvious replacements should be made
(``group" should be replaced with ``algebra" and so on).

\medskip

Similarly, we can prove an analogue of the small-subgroups theorem but we
restrict ourselves to an analogue of the finite-subgroup theorem.

\proclaim{Finite-dimensional-ideal theorem}.
If an algebra $G$ has a finite-dimensional two-sided ideal such that
the quotient algebra satisfies a given universal positive
closed first-order formula (in the language of algebras over the given
field), then $G$ has a characteristic
finite-dimensional two-sided ideal with the same property.

%%%%%%%%%%%%%%%%%%%%%%%%%%%%%%%%%%%%%%%%%%%%%%%%%%%%%%%%%%%%%%%%%%%%%%%%%%
\s 5.
The lattice of finite subgraphs

The word \emph{graph} in this section can be understood in any
reasonable sense:
all propositions are valid for directed and undirected graphs;
multiple edges and loops may be allowed or prohibited;
the vertices
and/or edges may be coloured. In the forbidden-subgraph theorem and
the local-embeddability theorem,
the word ``graph" may be even understood as ``hypergraph". All these
variations do not affect the proofs; of course, automorphisms
of a graph should be understood in the corresponding sense.

\proclaim{Forbidden-subgraph theorem}.
Let $\{\Gamma_1,\dots,\Gamma_l\}$ be a finite set of finite
graphs called \emph{forbidden} and considered up to isomorphism,
and let $G$ be some graph.
If $G$ contains a
 finite set $\=N$ of edges
such that $G\setminus \=N$
does not contain forbidden subgraphs, then
$G$ contains a finite
set of edges $\=H$
which is invariant with respect to all automorphisms
of~$G$  and has the same property: $G\setminus \=H$ does not
contain forbidden subgraphs. Moreover,
$
|\=H|\leq f^{t-1}(|\=N|),
$
where $f^k(x)$ is the $k$th iteration of the function $f(x)=x(x+1)$,
and $t$ is the maximal (in $i$) number of edges of $\Gamma_i$.
%\newline
\quad
In addition, if $\=H\ne\emptyset$, then $\=H\cap\=N\ne\emptyset$.

\Proof
It suffices to apply the main theorem to the lattice consisting of
cofinite subsets of the edge set of $G$. Clearly, this
lattice is Noetherian and the function
$$
\cod (X)\:=\hbox{the number of edges of the graph $G\setminus X$}
$$
satisfies all conditions from the definition of codimension.
The semigroup $\Phi$ is taken to be the automorphism group of~$G$ and
$\PP(N_1,\dots, N_t)$ is the following property:

{\sl\noindent
\narrower\narrower\narrower
$G$ contain no forbidden subgraph whose first edge
lies in $N_1$, second edge lies in $N_2$,\dots\

}

We assume that the edges of each
forbidden graph are enumerated somehow. Clearly, the property $\PP$ is
monotone:
$$
\PP(N_1,\dots, N_t) \imp \PP(N_1',\dots, N_t'),
\qbox{if
$N_i'\subseteq N_i$}.
$$
The multilinearity is also obvious:
$$
\big(
\PP(N_1',N_2\dots, N_t)\wedge \PP(N_1'',N_2\dots, N_t)
\big)
\imp
\PP(N_1'\cup N_1'',N_2\dots, N_t).
$$
It remains to apply the main theorem and note that
the property $\PP(N,\dots,N)$ means precisely the absence of forbidden
subgraphs in $N$.

Now, $\=H=G\setminus H$ and $\=N=G\setminus N$ intersect, because,
according to the main theorem, $H$ is contained in the sublattice
generated by all images of
$N$ under automorphisms of $G$.
In particular, ${H\supseteq\bigcap\limits_{\phi\in\Aut G}\phi(N)}$, i.e.
${\=H\subseteq\bigcup\limits_{\phi\in\Aut G}\phi(\=N)}$.
Therefore, if $\=H\ne\emptyset$, then ${\=H\cap\phi(\=N)\ne\emptyset}$
for some $\phi\in\Aut G$. By virtue of the invariance of $\=H$, this means
that $\=H\cap\=N\ne\emptyset$ as required.
The theorem is proven.

\medskip

This theorem can be significantly strengthened if we do not care about
the estimate. We say that a graph $X$
\emph{locally embeds} into a graph $Y$ if any finite subgraph of $X$
is isomorphic to some subgraph of $Y$.

\proclaim{Local-embeddability theorem}.
For any graph $G$ and any finite set $\=N$ of its edges,
there exists a finite
set of edges~$\=H$ invariant with respect to all automorphisms
of $G$ and such that the graph $H=G\setminus\=H$ locally embeds into
the
graph $N=G\setminus\=N$.

\Proof
Let $\Gamma_1,\Gamma_2,\dots$ be all finite graphs not embeddable into $N$.
We have to show that

{\sl\noindent\narrower
all subgraphs of $G$ isomorphic to $\Gamma_i$ can be destroyed
by removing a finite $(\Aut G)$-invariant set of edges $\=H$
provided we know that these subgraphs can be
destroyed
by removing some finite
set of edges $\=M$.

}

\noindent
To prove this assertion we use the induction on $|\=M|$. Let us start with
$\=M=\=N$.
By the forbidden-subgraph theorem, for each positive integer $n$,
there exists a finite set $\=H_n$ of edges such that
\item{1)}
$\=H_n$ is invariant with respect to all automorphisms of $G$;
\item{2)}
the
graphs $\Gamma_1,\dots,\Gamma_n$ are
not embeddable into $H_n=G\setminus\=H_n$;
\item{3)}
$\=H_n\cap\=M\ne\emptyset$ if $\=H_n\ne\emptyset$.

\enditem
If all sets $\=H_n$ are empty, then we have nothing to prove.
If there is a nonempty set $\=H_k$, then we consider the graph
$G'=H_k=G\setminus\=H_k$. This graph contains a
finite set
$\=M'=\=M\setminus(\=M\cap\=H_k)$ of edges such that
none of $\Gamma_i$ embeds into~$G'\setminus\=M'$.

Moreover $|M'|<|M|$ by property 3) of $\=H_k$. Therefore, by
the induction hypothesis, $G'$ contains a finite invariant set
$\=H'$ of edges such that none of $\Gamma_i$
embeds into $G'\setminus\=H'=G\setminus(\=H_k\cup\=H')$; this is
what we want, because $\=H_k\cup\=H'$ is an invariant set. Indeed,
$\=H_k$ is invariant with respect to all automorphisms of
$G$ by definition; $\=H'$ is invariant with respect to $\Aut G'$, and,
hence, with respect to $\Aut G$, because $G'$ is
an $(\Aut G)$-invariant subgraph of $G$. This completes the proof.

\Example.
Consider the undirected graph $G$ homeomorphic to the straight line and
the graph $N$ obtained from $G$ by removing one edge. Clearly,
it is impossible to remove a finite automorphism-invariant set of edges 
from $G$ in such a way that obtained graph 
$H$ is embeddable in $N$ (because the automorphism group of~$G$ acts 
transitively on edges). This example shows that
we cannot replace local embeddability with embeddability
in the theorem.

%\begin{sled}
\proclaim{Planarity theorem}.
If a graph can be made planar by removing a finite number of edges, then
it can be made planar by removing a finite set of edges
which is
invariant with respect to all automorphisms of the graph.

\Proof
By the Kuratowski--Erd\H os--Wagner theorem [Wag67]
a graph $G$ is planar if and only if
\-
the number of its edges is at most continuum;
\-
the number of its vertices of
degree larger than two is countable (or finite);
\-
it does not contain subgraphs homeomorphic to the complete graph on five
vertices $K_5$ or the complete bipartite graph with three vertices in each
fraction $K_{3,3}$ (Fig.~1); i.e.
$G$ does not contain
subgraphs isomorphic to graphs obtained from $K_5$
or $K_{3,3}$ by subdivisions of edges.

\enditem

\goodbreak
%\vskip1cm plus 1cm minus5mm
\bigskip
\centerline{\input 1.PIC}
\nobreak%
%\vskip5mm%
\centerline{Fig. \lowercase{1}}%
%\vskip1cm plus 1cm minus 5mm%
\goodbreak
\bigskip
%\par
%\noindent

The first two properties are not affected by
removing or adding a finite number of edges; the third property
is inherited by
locally embeddable graphs. Therefore, the assertion follows immediately
from the local-embeddability theorem.

\medskip

The following proposition shows that, in the planarity theorem, no
estimate of the cardinality of the invariant removed set of edges
is possible.

\Proposition.
For each positive integer $n$, there exists a finite graph $G_n$ which
becomes planar after removing five edges but cannot be
made planar by removing an invariant set consisting of less
than $n$ edges.

\Proof
Let us take the graph $K_5$ and subdivide each edge of a length-five
cycle
onto~$n$ parts. Now, let us glue together $n$ copies
of the obtained graph along the cycle of length $5n$
with rotations (Fig. 2).

\vfil\break

\goodbreak
%\vskip1cm plus 1cm minus5mm
\bigskip
\centerline{\input 2.PIC}
\nobreak%
%\vskip5mm%
\centerline{Fig. \lowercase{2}}%
%\vskip1cm plus 1cm minus 5mm%
\goodbreak
\bigskip
%\par
%\noindent

The obtained graph $G_n$ becomes planar after removing five
edges --- each $n$th edge on the $5n$-cycle (Fig. 3).
However, removing small invariant set of edges cannot make this
graph planar, because, among automorphisms of $G_n$, there is the rotation
through
one edge along the cycle of length $5n$ and, therefore, the orbit of each
edge has at least~$n$ elements. Thus, any invariant set
of edges containing less than $n$ elements must be empty. It remains to
note, that the graph~$G_n$ itself is not planar, because it contains
a subgraph homeomorphic to $K_5$.

%\vfil\break

\goodbreak
%\vskip1cm plus 1cm minus5mm
\bigskip
\centerline{\input 3.PIC}
\nobreak%
%\vskip5mm%
\centerline{Fig. \lowercase{3}}%
%\vskip1cm plus 1cm minus 5mm%
\goodbreak
\bigskip
%\par
%\noindent

In conclusion, we note that, at least for countable graphs, there is 
an analogue of the planarity theorem in which planarity is replaced
with embeddability into any fixed surface. To show this, it suffices to
recall a theorem of Erd\H os which says that a countable graph embeds
into a surface $S$ if and only if each its finite subgraph embeds into $S$.

%%%%%%%%%%%%%%%%%%%%%%%%%%%%%%%%%%%%%%%%%%%%%%%%%%%%%%%%%%%%%%%%%%%%%%%%%%
%\s 6.
%Другие приложения

%

%%%%%%%%%%%%%%%%%%%%%%%%%%%%%%%%%%%%%%%%%%%%%%%%%%%%%%%%%%%%%%%%%%%%%%%%%%
%\s 5.
%Топологические группы и группы Ли

%

%%%%%%%%%%%%%%%%%%%%%%%%%%%%%%%%%%%%%%%%%%%%%%%%%%%%%%%%%%%%%%%%%%%%%%%%%%
\vfil\break

\s 6.
Elementary mathematics

\Problem 1.
In the three-dimensional Euclidean space, there is a set $X$.
It is known that we can remove a finite set of
points from $X$ in
 such a way that no \the\year\ of the remaining points lie
on the same sphere. Show that this finite set can be chosen
invariant under all symmetries(= isometries)
of $X$.

\Solution
It suffices to apply the main theorem to
the lattice of cofinite subsets
of $X$ (this lattice is obviously Noetherian)
and take
$\Phi$ to be the symmetry group of $X$
and $\PP$ to be the following \the\year-linear
$\Phi$-invariant
monotone predicate:
$$
\PP(N_1,\dots,N_{\the\year})=
\hbox{(no points $x_1\in N_1,\dots,x_{\the\year}\in N_{\the\year}$
lie on the same sphere)}.
$$
On the lattice, there is also a natural codimension:
$\cod(X\setminus K)\:=|K|$
that makes it possible to estimate the cardinality of the symmetric removed set via
the cardinality of the initial (nonsymmetric) finite set.

\Problem 2.
There were chosen $10^{100}$ excellent
candidates
for a mission to Mars. The only problem is that they do not have enough
respect for each other. The organisers noted that expelling some ten
persons would result in an \emph{efficient team} in the sense that, among
any five of the remaining candidates, there is at least one respected by
the majority (of this five). Show that a nonempty efficient team can be
build \emph{fairly}, i.e. in such a way that the expelled set is invariant
under all permutations of candidates preserving the relation ``respects".
(Certainly, the binary relation ``respects" may be non-transitive,
non-symmetric, and even non-reflexive.)

\Solution
It suffices to apply the main theorem to the
lattice of all subsets of the
set of candidates $X$
putting $\Phi$ to be the group of permutations of $X$
preserving
the relation ``respects"
and $\PP$ to be the following pentalinear
$\Phi$-invariant
monotone predicate:
$$
\PP(N_1,\dots,N_5)=
\hbox{(any candidates $x_1\in N_1,\dots,x_5\in N_5$
form an efficient five)}.
$$
There is a natural codimension: $\cod(X\setminus K)\:=|K|$,
that makes it possible to estimate the number of fairly expelled candidates:
$$
\cod H \leq f^{t-1}(\cod N)=f^{5-1}(10)<
{11\over10}\cdot\({11\over10}\cdot\({11\over10}\cdot\({11\over10}\cdot
10^2\)^2\)^2\)^2=
\({11\over10}\)^{15}\cdot10^{16}=11^{15}\cdot10\ll10^{100},
$$
i.e. the remaining efficient team is nonempty.

%%%%%%%%%%%%%%%%%%%%%%%%%%%%%%%%%%%%%%%%%%%%%%%%%%%%%%%%%%%%%%%%%%%%%%%%%%%
\REFERENCES
\itemwidth1.8cm
\parindent2cm

\[BeK03]
Belyaev V.V., Kuzucuo\u glu M.
Locally finite barely transitive group
{//Algebra i Logika.} 2003. V.42.  no.3. P.261--270.

\[KaM82]
Kargapolov M.I., Merzlyakov Yu.I.
Fundamentals of group theory.
Moscow: ``Nauka", 1982.

\[KlM09]
Klyachko Ant. A., Melnikova Yu.B.
A short proof of the Khukhro--Makarenko theorem on large characteristic
subgroups with laws
{// Mat. Sbornik}, 2009, 200:5, 33--36.
See also
arXiv:0805.2747 .

\[Kur62]
Kurosh A.G.
Lectures on general algebra.
Moscow: ``Fiz.-Mat.Lit.", 1962.

\[She78]
Shemetkov L.A.
Formations of finite groups.
Moscow: ``Nauka", 1978.

\[AST13]
A. Arikan, H. Smith, N. Trabelsi,
On certain application of the Khukhro--Makarenko theorem,
Glasgow Math. J. 55(2013), 275--283.

\[BrNa04]
Bruno B., Napolitani F.
A note on nilpotent-by-\v Cernikov groups
{// Glasgow Math. J.} 2004. 46, 211-215.

\[Corn13]
Y. Cornulier
(http://mathoverflow.net/users/14094/yves-cornulier),
Large abelian characteristic subgroups in abelian-by-countable groups,
URL (version: 2013-12-15): http://mathoverflow.net/q/151889

\[Higg56]
P. J. Higgins,
Groups with multiple operators,
Proc. London Math. Soc. (3) 6 (1956), 366-416.

\[KhM07a]
Khukhro E.I., Makarenko N.Yu.
Large characteristic subgroups satisfying multilinear commutator identities
{// J. London Math. Soc.} 2007. {V.75}. no.3, P.635--646.

\[KhM07b]
Khukhro E.I., Makarenko N.Yu.
Characteristic nilpotent subgroups of bounded co-rank and
automorphically-invariant ideals of bounded codimension in Lie algebras
{// Quart. J. Math.} 2007. {V.58}. P.229--247.

\[KhM08]
Khukhro E.I., Makarenko N.Yu.
Automorphically-invariant ideals satisfying multilinear identities,
and group-theoretic applications
{// J. Algebra} 2008. {V.320}. no.4. P.1723--1740.

\[KhKMM09]
E. I. Khukhro, Ant. A. Klyachko, N. Yu. Makarenko, and Yu. B. Melnikova
Automorphism invariance and identities.
Bull. London Math. Soc. (2009), 41(5): 804-816.
See also arXiv:0812.1359

\[MSh12]
N.Yu. Makarenko, P. Shumyatsky,
Characteristic subgroups in locally finite groups,
Journal of Algebra, (2012), 352:1, 354--360.

\[Wag67]
K. Wagner,
Fastpl\"attbare Graphen,
J. Combinatorial Theory. (1967), 3, 326--365.

\end